\documentclass[12pt]{amsart}
\usepackage{amsmath,amscd,amssymb,amsfonts}
\newcommand{\ssb}{\raise.15ex\h{${\scriptscriptstyle\bullet}$}}
\newcommand{\ssc}{\,\raise.15ex\h{${\scriptstyle\circ}$}\,}
\newcommand{\D}{{\mathcal D}}
\newcommand{\M}{{\mathcal M}}
\newcommand{\Hc}{{\mathcal H}}
\newcommand{\OO}{{\mathcal O}}

\newcommand{\C}{{\mathbf C}}
\newcommand{\Q}{{\mathbf Q}}
\newcommand{\R}{{\mathbf R}}
\newcommand{\Z}{{\mathbf Z}}
\newcommand{\h}{\hbox}
\newcommand{\q}{\quad}

\newcommand{\nin}{\par\noindent}
\newcommand{\bs}{\par\bigskip}
\newcommand{\ms}{\par\medskip}
\newcommand{\sk}{\par\smallskip}
\newcommand{\msum}{\h{$\sum$}}

\newcommand{\Gr}{{\rm Gr}}

\newcommand{\MHM}{{\rm MHM}}
\newcommand{\MHW}{{\rm MHW}}
\newcommand{\into}{\hookrightarrow}
\newcommand{\simto}{\buildrel{\sim}\over\longrightarrow}

\newcommand{\les}{\leqslant}

\begin{document}
\title{On the definition of mixed Hodge modules}
\author[M. Saito]{Morihiko Saito}
\address{RIMS Kyoto University, Kyoto 606-8502 Japan}
\begin{abstract}
We give some details of a simpler definition of mixed Hodge modules which has been announced in some papers.
Compared with earlier arguments, this new definition is simplified by using Beilinson's maximal extension together with stability by subquotients systematically.
\end{abstract}
\maketitle
\centerline{\bf Introduction}
\bs\nin
Let $X$ be a complex algebraic variety (which is not assumed irreducible), and $A$ be a subfield of $\R$. Let $\MHW(X,A)$ be the abelian category of weakly mixed Hodge modules. More precisely, it consists of filtered regular holonomic $\D_X$-modules endowed with $A$-structure and the weight filtration
$$\M=((M;F,W),(K,W);\alpha)\in MF_{rh}W(\D_X,A),$$
such that the $\Gr_k^W\M$ are pure Hodge modules of weight $k$ (assumed always polarizable) for any $k\in\Z$.
(See (1.1) for $MF_{rh}W(\D_X,A)$.)
This means that a weakly mixed Hodge module is obtained by successive extensions of pure Hodge modules without imposing any conditions on the
extensions. (So it is called `weakly mixed'.)
Here $F$ and $W$ are respectively called the {\it Hodge} and {\it weight} filtrations. Note that any morphism of $\MHW(X,A)$ is bistrictly compatible with $F,W$ (see [Sa1, Prop.~5.1.14]).
\sk
The category of mixed Hodge modules $\MHM(X,A)$ is the abelian full subcategory of $\MHW(X,A)$ defined by increasing induction on the dimension of support as follows:
\sk
For $\M\in\MHW(X,A)$ with ${\rm supp}\,\M=X$, $\M$ belongs to $\MHM(X,A)$ if and only if, for any $x\in X$, there is a Zariski-open neighborhood $U_x$ of $x$ in $X$ together with a function $g$ on $U_x$ such that $U_x':=U_x\setminus g^{-1}(0)$ is smooth and dense in $U_x$, and the following two conditions are satisfied:
\ms\nin
(C1)\, The restriction $\M':=\M|_{U_x'}$ is an admissible variation of mixed $A$-Hodge structure in the sense of [Ka2], [SZ] (up to a shift of complex).
\ms\nin
(C2)\, The nearby and vanishing cycle functors along $g=0$ are well-defined for $\M|_{U_x}$ (see (2.2)), and $\varphi_{g,1}\M|_{U_x}\in\MHM(g^{-1}(0),A)$ (where $\dim g^{-1}(0)<\dim X$ by assumption).
\ms
Note that admissible variations of mixed Hodge structures are mixed Hodge modules under this definition since $g^{-1}(0)$ may be empty in condition (C2).
However, we need the following theorems in order to justify this.
\bs\nin
{\bf Theorem 1.} {\it Conditions {\rm (C1-2)} is independent of the choice of $U_x$, $g$.
More precisely, assume they are satisfied for some $U_x$, $g$ for any $x$.
Then {\rm (C2)} is satisfied for any $U_x$, $g$, and {\rm (C1)} is satisfied in case the underlying perverse sheaf of $\M'$ is a local system.}
\bs\nin
{\bf Theorem 2.} {\it The categories $\MHM(X,A)$ for complex algebraic varieties $X$ are stable by the canonically defined cohomological functors $\Hc^jf_*$, $\Hc^jf_!$, $\Hc^jf^*$, $\Hc^jf^!$, $\psi_g$, $\varphi_{g,1}$, $\boxtimes$ for morphisms of complex algebraic varieties $f$ and functions $g$ on complex algebraic varieties. Moreover these functors are compatible with the corresponding functors of the underlying perverse sheaves with $A$-coefficients.}
\bs
Compared with earlier arguments, this new argument is simplified by using Beilinson's maximal extension together with the stability by subquotients systematically.
Combining Theorem~2 with the construction in [Sa3], we get the following.
\ms\nin
{\bf Corollary 1.} {\it There are canonically defined functors $f_*$, $f_!$, $f^*$, $f^!$, $\psi_g$, $\varphi_{g,1}$, $\boxtimes$, $\otimes$, ${\Hc}om$ between the bounded derived categories $D^b\MHM(X,A)$ for complex algebraic varieties $X$ so that we have the canonical isomorphisms $H^jf_*=\Hc^jf_*$, etc., where $f$ is a morphism of complex algebraic varieties, $g$ is a function on a complex algebraic variety, $H^j$ is the usual cohomology functor of the derived categories, and $\Hc^jf_*$, etc.\ are as in Theorem~$2$. Moreover the above functors between the $D^b\MHM(X,A)$ are compatible with the corresponding functors of the underlying complexes with $A$-coefficients.}
\ms
This work is partially supported by Kakenhi 24540039.
\bs\bs
\centerline{\bf 1. Preliminaries}
\bs\nin
{\bf 1.1. Filtered regular holonomic $\D$-modules with $A$-structure.} Let $X$ be a complex algebraic variety, and $A$ be a subfield of $\R$. An object of the category $MF_{rh}W(\D_X,A)$ consists of
$$((M;F,W),(K,W);\alpha),$$
where $(M,F)$ is a filtered regular holonomic $\D_X$-module with a finite filtration $W$, $(K,W)$ is a filtered perverse sheaf with coefficients in $A$, and $\alpha$ is an isomorphism between ${\rm DR}(M)$ and $K\otimes_A\C$ compatible with $W$.
\ms\nin
{\bf 1.2.~Filtered $\D$-modules on singular varieties.} A filtered regular holonomic $\D_X$-module $(M,F)$ on a singular complex algebraic variety $X$ is defined by taking closed embeddings into smooth varieties $U\into Y$, and using filtered regular holonomic $\D_Y$-modules $(M_Y,F)$ where $U$ are open subvarieties of $X$. 
We use right $\D$-modules here since it behaves better than left $\D$-modules under the direct images by closed embeddings (see (1.3) below). We assume
$$\h{The $\Gr^F_pM_Y$ are annihilated by the ideal of the support of $M_Y$.}
\leqno(1.2.1)$$
This condition is satisfied in the case of the underlying filtered $\D$-modules of mixed Hodge modules by [Sa1, Lemma~3.2.6].
\sk
If we have two closed embeddings $U\into Y_a$ ($a=1,2)$ with $Y_a$ smooth and if there is a smooth morphism $\pi:Y_1\to Y_2$ compatible with the closed embeddings (by replacing $Y_1$ with $Y_1\times Y_2$), then we have a relation between the two filtered regular holonomic $\D_{Y_a}$-modules $(M_{Y_a},F)$ ($a=1,2$) by a given isomorphism
$$\pi_*(M_{Y_1},F)\cong(M_{Y_2},F),$$
(satisfying some compatibility condition, see e.g. [Sa1, 2.1.20]).
This implies that $(M_{Y_2},F)$ is determined by $(M_{Y_1},F)$.
Here $\pi_*(M_{Y_1},F)$ is the direct image as a filtered $\D$-module, and it is a filtered $\D$-module by using a local section $\sigma$ of $\pi$ compatible with the closed embeddings (by using the minimal embedding), since $\pi_*\ssc\sigma_*=id$. This implies also
$$(M_{Y_1},F)\cong\sigma_*(M_{Y_2},F),$$
which shows that $(M_{Y_1},F)$ is determined by $(M_{Y_2},F)$.
We can also show that $\sigma_*(M_{Y_2},F)$ is independent of the choice of the local section $\sigma$ by (1.2.1) and [Sa1, Lemma~3.2.6].
\ms\nin
{\bf 1.3.~Direct images of filtered $\D$-modules.} In case of a closed embedding of smooth complex algebraic varieties $f:X\to Y$, the direct image of a filtered $\D$-module is locally given by the tensor over $\D_X$ with
$$\D_{X\to Y}:=\OO_X\otimes_{f^{-1}\OO_Y}f^{-1}\D_Y,$$
for right $\D$-modules. (If we use left $\D$-modules, we need the the twist by $\omega_{X,Y}$ together with a shift of the filtration $F$.)
If $X$ is locally defined by $y_1=\cdots=y_r$ with $y_1,\dots,y_m$ local coordinates of $Y$, it is locally given by the tensor with $\C[\partial_1,\dots,\partial_r]$ over $\C$, where $\partial_i:=\partial/\partial y_i$.
\sk
In case of a smooth projection, the direct image is defined by using the relative de Rham complex ${\rm DR}_{X/Y}$ (up to a shift of complex) which is locally the Koszul complex associated with the action of $\partial/\partial x_i$ if $x_1,\dots,x_n$ are local coordinates of the fiber.
(Note, however, that this does not work for a smooth morphism which is not a smooth {\it projection} unless we consider cohomological direct image {\it sheaves}, instead of direct image {\it complexes}.)
\sk
In general, the direct image of filtered $\D$-modules by a morphism of smooth varieties $f:X\to Y$ is defined by using the factorization $f=pr_2\ssc i_f$, where $i_f:X\to X\times Y$ is the graph embedding, and $pr_2$ is the second projection.
(It may be simpler to use this construction instead of induced $D$-modules as in [Sa1, 3.3.6] also for the definition of the direct image of the $V$-filtration.)
\sk
It is easy to see that the above direct image in the general case is naturally isomorphic to the complex of the induced $\D_Y$-module associated with the (sheaf-theoretic) direct image of the filtered differential complex ${\rm DR}_X(M,F)$. This implies the compatibility between the direct images of filtered differential complexes and filtered $\D$-modules.
\sk
The direct image for a morphism of singular varieties is rather complicated.
For the direct image of mixed Hodge modules, we may assume that the morphism is projective by using a Beilinson-type resolution, and the cohomological direct image is actually enough. So it is reduced to the case of a morphism of smooth varieties.
\ms\nin
{\bf 1.4.~Left and right $\D$-modules.}
The transformation between filtered left and right $\D_Y$-modules on a complex manifold $Y$ of dimension $m$ is given by associating to a filtered $\D_Y$-module $(M,F)$
$$(\Omega_Y^m,F)\otimes_{\OO_Y}(M,F),
\leqno(1.4.1)$$
where the filtration $F$ on $\Omega_Y^m$ is defined by the condition $\Gr^F_p\Omega_Y^m=0$ for $p\ne -m$. By choosing local coordinates $y_1,\dots,y_m$, the sheaf $\Omega_Y^m$ is locally trivialized by $dy_1\wedge\cdots\wedge dy_m$, and the transformation forgetting $F$ is given by the anti-involution $^*$ of $\D_Y$ defined by the conditions:
$$(\partial/\partial y_i)^*=-\partial/\partial y_i,\q g^*=g\,\,\,(g\in\OO_Y),\q(PQ)^*=Q^*P^*\,\,\,(P,Q\in\D_Y).
\leqno(1.4.2)$$
\ms\nin
{\bf 1.5.~Filtered differential complexes.} Let $Y$ be a smooth complex algebraic variety or a complex manifold. There is an equivalence of categories (see [Sa1, 2.2.10])
$${\rm DR}^{-1}_Y:D^bF(\OO_Y,{\rm Diff})\simto D^bF(\D_Y),$$
with quasi-inverse given by the de Rham functor $\rm DR$.
Here the source is the derived category of bounded filtered differential complexes, and the target is the bounded derived category of filtered $\D_Y$-modules $D^bF(\D_Y)$. Its proof is quite elementary, and follows from a rather trivial calculation of certain Koszul complexes.
\sk
In fact, for filtered $\OO_Y$-modules $(L,F),(L',F)$ with $F$ bounded below, $\xi:(L,F)\to(L',F)$ is called a filtered differential morphism if the following composition is a differential morphism of order $\les q$ for any $p,q$:
$$F_pL\to L\buildrel{\xi}\over\to L'\to L'/F_{p-q-1}L'.$$
It is easy to see that this condition is equivalent to that the corresponding morphism of induced $\D$-modules
$${\rm DR}_Y^{-1}\xi:(L,F)\otimes_{\OO_Y}(\D_Y,F)\to(L',F)\otimes_{\OO_Y}(\D_Y,F)$$
is a morphism of filtered right $\D$-modules, where
$${\rm DR}_Y^{-1}(L,F):=(L,F)\otimes_{\OO_Y}(\D_Y,F).$$
These are filtered right $\D_Y$-module, and are called induced modules.
\ms\nin
{\bf 1.6.~Derived categories.} The bounded derived category of filtered $\D$-modules $D^bF(\D_Y)$ are defined by using Verdier's theory [Ve1] (i.e. by inverting filtered quasi-isomorphisms, or equivalently, dividing by the subcategory consisting of filtered acyclic complexes). The category $D^bF(\OO_Y,{\rm Diff})$ can be defined similarly.
\sk
We needed the category of graded $\bigoplus_pF_p\D_Y$-modules in [Sa1] in order to apply Deligne's theory of decompositions associated to certain actions satisfying the hard Lefschetz property. Here we get a pair of decompositions for the direct images of the underlying filtered $\D$-module and perverse sheaf, and their compatibility follows from Deligne's theorem about a unique choice of decompositions. However, this does not imply a decomposition in the derived category of mixed Hodge modules.
\sk
This point was improved in [Sa3, 4.5.4] where we showed the decomposition theorem in the derived category of mixed Hodge modules by using a property of pure complexes similar to the $l$-adic case, and the decomposition theorem holds under the assumption that the morphism is proper, instead of projective.
\ms\nin
{\bf 1.7.~Dual.} In order to define a polarization as in [Sa1], we need that the dual ${\bf D}(M,F)$ of the underlying filtered $\D$-module $(M,F)$ is a filtered $\D$-module. So we have to show that $(M,F)$ is {\it Cohen-Macaulay}, i.e. $\Gr^FM$ is Cohen-Macaulay over $\Gr^F\D_Y$, see [Sa1, 5.1.13].
This is used in an essential way for the inductive argument in the proof of the stability by projective direct images when we consider the induced duality on the nearby and vanishing cycles.
\ms\nin
{\bf 1.8.~Remark.}  It is rather nontrivial to calculate the nearby and vanishing cycle functors for filtered $\D$-modules even in the normal crossing case [Sa3], and to prove the compatibilities of the induced dualities on the nearby and vanishing cycle functors for filtered $\D$-modules and for perverse sheaves [Sa2].
This does not seem to be necessarily easier than some problems of the associated mixed Hodge structures. (Here note that Kashiwara's original proof of the assertion in the Appendix of [Sa3] used a reduction to the two dimensional case, and some simplification was made by the writer of the Appendix. Note also that the position of the Weil operator $C$ in the last formula of [Sa3, A4] is different from [Ka2, 1.3.2]. In fact, these respectively follow the convention of Deligne, Steenbrink and that of Cattani, Kaplan. They produce the difference between $i^{q-p}$ and $i^{p-q}$ when we consider the associated Hermitian form for $u=v$ with type $(p,q)$.)
The above results are needed in order to show that polarizable variations of Hodge structure on smooth open subvarieties are uniquely extended to pure Hodge modules (see [Sa3, Th.~3.21]) and admissible variations of mixed Hodge structure are mixed Hodge modules.
(Note that pure Hodge modules are mixed Hodge modules since polarizable variations of Hodge structure are admissible variations of mixed Hodge structures.)
\bs\bs
\centerline{\bf 2. Extensions of weakly mixed Hodge modules}
\bs\nin
{\bf 2.1.~Pure Hodge modules.} For an irreducible variety $X$ of dimension $d$, we have an equivalence of categories
$${\rm MH}_X(X,n,A)={\rm VHS}_{\rm gen}(X,n-d,A)^p,
\leqno(2.1.1)$$
where the left-hand side is the category of polarizable Hodge modules of weight $n$ with strict support $X$ (i.e. the underlying perverse sheaf is the intersection complex with local system coefficients), and the right-hand side is the category of polarizable variations of Hodge structure of weight $n-d$ defined on smooth open subvarieties of $X$.
\sk
For a pure Hodge modules with strict support $X$, the Hodge filtration $F$ on the underlying $\D_Y$-module $M$ (defined for each local embedding $X\into Y$) is uniquely determined by its restriction to any dense open subvariety of $X$ by using [Sa1, Prop.~3.2.2(ii)] which implies
$$F_p(i_g)_*M=\msum_{i\ge 0\,}\bigl(V^{>-1}(i_g)_*M\cap j_*j^{-1}F_{p-i}(i_g)_*M\bigr)\partial_t^i,
\leqno(2.1.2)$$
where $g$ is any locally defined nonzero function, $(i_g)_*M$ is the direct image as a filtered $\D$-module by the graph embedding $i_g$, $t$ is the coordinate of $\C$, and $V$ is the filtration of Kashiwara [Ka1] and Malgrange [Ma].
In this paper, the filtration $V$ is indexed by $\Q$ so that the action of $\partial_tt+\alpha$ on $\Gr_V^{\alpha}$ is nilpotent for right $\D$-modules (i.e. $t\partial_t-\alpha$ is nilpotent for left $\D$-modules, see (1.4.2)).
\sk
In loc.~cit., the following conditions (i.e. (3.2.2-3) in loc.~cit.) are used:
$$t:(V^{\alpha}(i_g)_*M,F)\simto(V^{\alpha+1}(i_g)_*M,F[-1])\q(\alpha>-1),
\leqno(2.1.3)$$
$$\partial_t:(\Gr_V^{\alpha}(i_g)_*M,F)\simto(\Gr_V^{\alpha-1}(i_g)_*M,F[-1])\q(\alpha<0).
\leqno(2.1.4)$$
Note that the strict surjectivity of (2.1.4) for $\alpha=0$ is also needed in order to show (2.1.2) by using Prop.~3.2.2(ii) in loc.~cit., but this is not included in conditions (2.1.3-4), and we have to prove that it underlies a morphism of $\MHW(X,A)$ in order to prove its strict surjectivity.
\ms\nin
{\bf 2.2.~Well-definedness of nearby and vanishing cycle functors.}
Let $g$ be a function on $X$.
We say that the nearby and vanishing cycle functors along $g=0$ are {\it well-defined} for $\M\in\MHW(X,A)$ (i.e. $\M$ is {\it specializable} along $g^{-1}(0)$) if the following two conditions are satisfied:
\ms\nin
(W1)\, The three filtrations $F,W,V$ on the direct image of the underlying filtered $\D$-module by the graph embedding $i_g$ are compatible filtrations.
\ms\nin
(W2)\, There is the relative monodromy filtration $W$ for the action of the nilpotent part $N$ of the monodromy on $\psi_gK[-1]$, $\varphi_{g,1}K[-1]$ with respect to $L:=\psi_gW[-1]$ and $\varphi_{g,1}W$.
\ms\nin
{\bf 2.3.~Open direct images.}
Let $D$ be a locally principal divisor on $X$. We say that the open direct images $j_!.j_*$ by the inclusion $j$ of the complement of $D$ is well-defined for $\M'\in\MHW(X\setminus D,A)$, if there are $\M'_!$, $\M'_*\in\MHW(X,A)$ whose underlying perverse sheaves are respectively isomorphic to $j_!K'$, $j_*K'$ (where $K'$ is the underlying perverse sheaf of $\M'$), and the following condition is satisfied:
\ms\nin
(W3)\, For any locally defined function $g$ such that $g^{-1}(0)_{\rm red}=D_{\rm red}$, the nearby and vanishing cycle functors along $g=0$ are well-defined for $\M'_!$, $\M'_*$.
\ms
If this condition is satisfied, then $\M'_!$ and $\M'_*$ are respectively denoted by $j_!\M'$ and $j_*\M'$.
If $\M'=j^{-1}\M$ with $\M\in\MHW(X,A)$ and the nearby and vanishing cycle functors along $g=0$ are well defined for $\M$, then we have the canonical morphisms (see [Sa3, Prop.~2.11])
$$j_!j^{-1}\M\to\M\to j_*j^{-1}\M.
\leqno(2.3.1)$$
\ms\nin
{\bf 2.4.~Remarks.}~(i) The Hodge and weight filtrations on the open direct images are uniquely determined by the condition that the nearby and vanishing cycle functors are well-defined along $g=0$ where $g$ is any locally defined function such that $g^{-1}(0)_{\rm red}=D_{\rm red}$. (Note that this might depend on the choice of $g$.)
The formula for $F$ is similar to (2.1.2), see [Sa1, 3.2.3.1-2].
However, the formula for $W$ is rather complicated, and requires a generalization of a formula of Steenbrink and Zucker [SZ] and its dual (see [Sa3], Cor.~1.9]).
This is a special case of the Verdier-type extension theorem [Ve2] for weakly mixed Hodge modules, see [Sa2, Prop.~2.8].
\sk
For the proof of Theorem~2, we have to show that the open direct images are well-defined in the above sense. This can be reduced to the normal crossing case, and follows from [Sa2, Prop.~3.26].
\ms
(ii) In condition (W2) the shift of the filtration in $\psi_gW[-1]$ comes from the fact (or the normalization) that for a pure Hodge module $\M$ of weight $n$, the weight filtration of $\psi_g\M$ and $\varphi_{g,1}\M$ are the monodromy filtrations with center $n-1$ and $n$ respectively. We have the difference in the center of symmetry, since $\varphi_{g,1}\M$ is identified with the coimage of
$$N:\psi_{g,1}\M\to\psi_{g,1}\M(-1)$$
by the canonical morphism
$${\rm can}:\psi_{g,1}\M\to\varphi_{g,1}\M,$$
(i.e. the kernels of the two morphisms coincide and the last morphism is surjective).
\sk
For example, if $K=A_X[n]$ with $X$ smooth of dimension $n$ and $g$ is a coordinate, then $\psi_gA_X[n]=A_Y[n-1]$ with $Y=g^{-1}(0)$. There is no shift of filtration for $\varphi_{g,1}$ since $\varphi_{g,1}\M=\M$ in case $\M$ is supported in $g^{-1}(0)$.
\sk
Note also that the canonical morphism `${\rm can}$' is induced by $\partial/\partial t$ in the notation of (2.1.4), and we shift the Hodge filtration $F$ on $\Gr^0_V(i_f)_*M$ by 1 in case of right $\D$-modules so that the morphism preserves $F$. (In case of left $\D$-module, we shift $F$ on $\Gr^{-1}_V(i_f)_*M$ by $-1$, and there is no shift for $\Gr^0_V(i_f)_*M$. Note that $F$ is shifted by 1 when we take the direct image by $i_g$ for left $\D$-modules.)
\ms
(iii) It is not quite clear whether the above shift of weights on the nearby cycles is made properly in some papers of the twistor theory where the compatibility of sign is also nontrivial. These are essential in the proof of the stability by the projective direct images. In fact, if there is no shift of the filtration in the definition, then the twistor $\D$-module corresponding to the constant sheaf $A_X$ on a smooth variety $X$ (by the shifted de Rham functor) has always weight 0 independently of $\dim X$, and we would have a problem even in the case of a mixed twistor $\D$-module corresponding to the constant sheaf on a normal crossing varieties, since it should be pure if one follows the above convention. Also the monodromy filtration on the nearby cycles of the constant sheaf in the reduced normal crossing case cannot coincide with the weight filtration up to a shift (unless the weight is changed under the direct images by closed embeddings).
Note that the constant sheaf on the singular fiber is a subobject of the nearby cycles in the normal crossing case. This implies that the coprimitive part of the graded pieces of the nearby cycles are constant sheaves on closed strata of a fixed dimension in the case of mixed Hodge modules, and this should also hold for twistor $\D$-modules.
\sk
In the general case, however, the direct factors of each coprimitive part may have strict supports of various dimensions. This becomes more complicated by taking its direct image by the projective morphism
(which is needed when we calculate the $E_1$-term of the weight spectral sequence).
Here the direct factor with a fixed strict support on the base $Y$ is a direct sum of the direct images of the direct factors with various strict supports on $X$. Even if we have a positivity on each direct factor, it is quite nontrivial whether the signs of these positivities coming from various strict supports on $X$ are really compatible, and give a positivity on the direct sum. This is one of the most difficult point in the proof of the decomposition theorem for Hodge modules.
\ms
(iv) In condition (W2), $L$ is for the `limit weight filtration', and the relative monodromy filtration gives the weight filtration $W$ of $\psi_g\M$, $\varphi_{g,1}\M$ in condition (C2). (So it is sometimes called the relative monodromy weight filtration in the literature.) Note also that the limit Hodge filtration gives the Hodge filtration on the nearby and vanishing cycles, contrary to the case of the weight filtration.
\bs\bs
\centerline{\bf 3. Normal crossing case}
\bs\nin
{\bf 3.1.~Mixed Hodge modules with normal crossing singularities.} Let $Y$ be a smooth complex algebraic variety, and $D$ be a divisor with simple normal crossings.
Let $\MHM(Y,A)^{\rm nc}$ be the full subcategory of $\MHW(Y,A)$ consisting of $\M$ such that the restriction of $\M$ to each stratum of the stratification of $X$ associated with $D$ is locally constant, and the following condition is satisfied:
\sk
For any $y\in{\rm supp}\,\M\cap D$, let $D_1,\dots,D_r$ be the irreducible components of $D$ passing through $y$, and $V_{(i)}$ be the $V$-filtration along $D_i$. Then
\ms\nin
(N1)\, The $r+2$ filtrations $F,W,V_{(1)},\dots,V_{(r)}$ on the underlying $\D_Y$-module $M$ are compatible filtrations.
\ms\nin
(N2)\, For any $\alpha_i\in\Q\cap[-1,0]$ ($i\in[1,r]$) and $I\subset[1,r]$, there is the relative monodromy filtration $W(I)$ for the action of $\sum_{i\in I}N_i$ on $\Gr_{V_{(1)}}^{\alpha_1}\cdots\Gr_{V_{(r)}}^{\alpha_r}(M,W)$, satisfying the condition: $N_iW(I)_k\subset W(I)_{k-2}$ for any $k\in\Z$ and $i\in I$, see also [Ka2].
\ms\nin
{\bf 3.2.~Remarks.} (i) In case $\M$ is a pure Hodge module, condition (N1) is always satisfied, and we have the inductive structure for the relative monodromy filtrations by [CK].
The inductive structure in the mixed case is shown in [Ka2, Prop.~5.2.5], and a converse is proved in loc.~cit., Th.~4.4.1.
These imply that $\MHM(Y,A)^{\rm nc}$ coincides with the full subcategory of $\MHM(Y,A)$ defined by condition (3.25.1) in [Sa3], and moreover it is stable by external products. 
\ms
(ii) Condition (N1) assures that $F,W$ on $\Gr_{V_{(1)}}^{\alpha_1}\cdots\Gr_{V_{(r)}}^{\alpha_r}M$ are independent of the order of the $\Gr_{V_{(i)}}^{\alpha_i}$ for $i\in[1,r]$. It corresponds to the freeness of $\Gr_F^p\Gr_k^WM$ in case $M$ is the Deligne extension of an admissible variation of mixed Hodge structures with unipotent local monodromies.
\ms
(iii) Admissible variations of mixed Hodge structures with unipotent local monodromies on curves are defined by the freeness of $\Gr_F^p\Gr_k^WM$ in the Deligne extension case in Remark~(ii) together with the existence of the relative monodromy filtration, see [SZ]. The general case is reduced to the curve case in [Ka2] by considering the pull-backs by any morphisms from curves (including ramified coverings).
\ms
(iv) The above freeness is not sufficient in the curve case if the local monodromies are non-unipotent unless ramified coverings are taken. It can be replaced by the compatibility of three filtrations $F,W,V$ since $V$ is essentially given by the $\widetilde{t}$-adic filtration on some ramified covering with coordinate $\widetilde{t}$.
\ms
(v) There is no equivalence of categories between the category consisting of germs of weakly mixed Hodge modules around $x$ satisfying conditions (N1-2) and the category consisting of certain diagrams of mixed Hodge structures satisfying some conditions as in [Ka2]. In fact, after fixing the coordinates $x_i$, the latter should correspond to the category of `monodromical' mixed Hodge modules on $\C^n$, where monodromical means that the Hodge filtration $F$ is compatible with the decomposition of the $\D$-module obtained by the eigenspaces of the commuting actions of $x_i\partial/\partial x_i$.
\sk
This implies that the compatibility of a decomposition with the Hodge filtration $F$ is rather nontrivial even in case we get a decomposition for the corresponding combinatorial data of a mixed Hodge module with normal crossing singularity.
In order to solve this problem, we need for instance [Sa1, Lemma~5.1.4].
\bs\bs
\centerline{\bf 2. Proof of the main theorems}
\bs\nin
{\bf 4.1.~Proposition.} {\it Conditions {\rm (W1-2)} in $(2.2)$ are stable by subquotients in $\MHW(U,A)$.}
\ms\nin
{\it Proof.} See [Sa3, Prop.~2.5].
\ms
This implies the following.
\ms\nin
{\bf 4.2.~Corollary.} {\it The full subcategory $\MHM(X,A)$ is stable by subquotients in $\MHW(X,A)$.}
\ms\nin
{\bf 4.3.~Proposition.} {\it Condition $(W1$-$2)$ in $(2.2)$ are stable by the cohomological direct images $\Hc^jf_*$ under projective morphisms $f$, and we have}
$$\Hc^jf_*\,\psi_{fg}=\psi_g\,\Hc^jf_*,\q\Hc^jf_*\,\varphi_{fg,1}=\varphi_{g,1}\,\Hc^jf_*.$$
\ms\nin
{\it Proof.} See [Sa3, Th.~2.14].
\ms\nin
{\bf 4.4.~Proposition.} {\it In the notation of $(3.1)$, let $g$ be a function on $Y$ such that $D':=g^{-1}(0)\subset D$. Let $j:Y\setminus D'\into Y$ be the inclusion. Then, for any $\M\in\MHM(Y,A)^{\rm nc}$, the nearby and vanishing cycle functors along $g=0$ are well-defined, the open direct images $j_!,j_*$ are well-defined for $j^{-1}\M$, and the full subcategory $\MHM(Y,A)^{\rm nc}\subset \MHW(Y,A)$ in $(3.1)$ is stable by the functors
$$\psi_g,\q\varphi_{g,1},\q j_!j^{-1},\q j_*j^{-1}.$$
Moreover, for any admissible variation of mixed Hodge module $\M'$ on $Y\setminus D'$ which is identified with a weak mixed Hodge module, $j_!$, $j_*$ are well defined for $\M'$, and $j_!\M'$, $j_*\M'$ belong to $\MHM(Y,A)^{\rm nc}$.}
\ms\nin
{\it Proof.} This follows from [Sa3, Prop.~3.25-26 and the proof of Th.~3.27].
\ms\nin
{\bf 4.5.~Proposition.} {\it In the notation of conditions {\rm (C1-2)} in the introduction, $\M|_U$ is a subquotient of the direct sum of
$$(j_{U',U})_!(\M'\otimes g^*E_{0,b}),\q\varphi_{g,1}\M|_U,
\leqno(4.5.1)$$
for $b\gg 0$, where $j_{U',U}:U'\into U$ is the inclusion, and $E_{a,b}$ is the variation of mixed $A$-Hodge structure of rank $b-a+1$ on $\C^*$ such that $\Gr^W_kE_{a,b}=A(-j)$ if $k=2j$ with $j\in\Z\cap[a,b]$, and $0$ otherwise.}
\ms\nin
{\it Proof.}
Note first that $\M'\otimes g^*E_{0,b}$ is an admissible variation of mixed Hodge structure on $U'$ (up to a shift of complex), and $(j_{U',U})_!(\M'\otimes g^*E_{0,b})$ is well-defined by Propositions~(4.3) and (4.4) together with a resolution of singularities.
\sk
As is known in Beilinson's argument for gluing of perverse sheaves [Be2], it is rather easy to show (see Lemma~(4.6) below) that $\M|_U$ is isomorphic to the cohomology of the complex (see also [Sa3], p.~267)
$$\psi_{g,1}\M|_U\to\Xi_g\M'\oplus\varphi_{g,1}\M|_U\to\psi_{g,1}\M|_U(-1),
\leqno(4.5.2)$$
where $\Xi_g\M'$ is Beilinson's maximal extension of $\M'$, and is defined by
$$\Xi_g:={\rm Ker}\bigl((j_{U',U})_!(\M'\otimes g^*E_{0,b})\to(j_{U',U})_*(\M'\otimes g^*E_{1,b})\bigr)\q(b\gg 0),$$
see [Sa5] and also [Sa4].
This implies that $\M|_U$ is a subquotient of the direct sum of
$$\Xi_g\M',\q\varphi_{g,1}\M|_U,$$
and the assertion follows from the above definition of $\Xi_g\M'$.
\ms\nin
{\bf 4.6.~Lemma.} {\it The cohomology of $(4.5.2)$ is isomorphic to $\M|_U$.}
\ms\nin
{\it Proof.} Let $\M$ and $j$ respectively denote $\M|_U$ and $j_{U',U}$ to simplify the notation. Then $\varphi_{g,1}\M$ is given by the inductive limit of the single complex associated with the double complex
$$\begin{CD}j_!\M'@>>>\M\\
@VVV @VVV\\
j_!(\M'\otimes g^*E_{0,b})@>>>j_*(\M'\otimes g^*E_{0,b})
\end{CD}$$
(This can be proved by applying the functor $\varphi_{g,1}$ to the complex, see also [Sa4, Prop.~1.5].)
\sk
Moreover, $\Xi_g\M'$, $\psi_{g,1}\M$, $\psi_{g,1}\M(-1)$ are respectively given by the inductive limit of the complexes
$$\aligned j_!(\M'\otimes g^*E_{0,b})\to j_*(\M'\otimes g^*E_{1,b}),\\
j_!(\M'\otimes g^*E_{0,b})\to j_*(\M'\otimes g^*E_{0,b}),\\
j_!(\M'\otimes g^*E_{1,b})\to j_*(\M'\otimes g^*E_{1,b}).\endaligned$$
Then it is easy to verify that (4.5.2) is quasi-isomorphic to $\M$ (up to a shift of complex).
This finishes the proof of Lemma~(4.6).
\ms\nin
{\bf 4.7.~Proposition.} {\it In the notation of $(3.1)$, we have}
$$\MHM(Y,A)^{\rm nc}\subset\MHM(Y,A).
\leqno(4.7.1)$$
\ms\nin
{\it Proof.}
This follows by induction on the dimension of the support of $\M$ and the number of the irreducible components with the maximal dimension by using Propositions~(4.1), (4.4) and (4.5).
\ms\nin
{\bf 4.8.~Proof of Theorem~1.} By Propositions~(4.1) and (4.5) together with induction on the dimension of the support, the assertion is reduced to the case where $\M$ is the open direct image of an admissible variation of mixed Hodge structure by the inclusion of the complement of a divisor.
Here the open direct image can be defined by using Propositions~(4.3) and (4.4) together with a resolution of singularities. Moreover, these propositions imply that the assertion is further reduced to the case where the union of the divisor and $g^{-1}(0)$ has normal crossings.
In fact, for a birational projective morphism of smooth complex algebraic varieties $f:X\to Y$ and for a variation of mixed Hodge structure $\M$ on $Y$ which is identified with a weakly mixed Hodge module on $Y$, we have a canonical injective morphism
$$\M\to\Hc^0f_*\M',$$
where $\M'$ is the pull-back of $\M$ to $X$.
(For the proof of this we cannot use the adjunction since this is not yet proved.)
For instance, let $j:U\into Y$ be an affine open subvariety of $Y$ over which $f$ is an isomorphism.
We have canonical morphisms
$$j_!j^{-1}\M\to \M,\q j_!j^{-1}\M\to\Hc^0f_*\M',$$
where $j_!j^{-1}\M$ is defined by using Proposition~(4.3) and (4.4). We also apply these to $\M'$, and take the direct image by $f$ to get the second morphism.
So it is enough to show that the kernel of the first morphism is contained in that of the second, since the first morphism is surjective.
But this can be reduced to the assertion for the underlying perverse sheaves, where the assertion is well-known.
\sk
We may then assume that the union of the two divisors have normal crossings by using Propositions~(4.1) and (4.3).
So the assertion follows from Proposition~(4.4).
\ms\nin
{\bf 4.9.~Proof of Theorem~2.}
By the construction of the cohomological functors in [Sa3] (and using Proposition~(4.3)), it is enough to show the stability by the following functors: the direct images under the inclusions of the complements of principal divisors, the nearby and vanishing cycle functors, and the external products. These are exact functors and compatible with the passages to the subquotients of (4.5.1). Hence we can apply these functors to (4.5.1) instead of $\M|_U$, and proceed by induction on the dimension of the support.
For the external product $\boxtimes$, we use induction on the dimension of the support of the external product together with the stability of $\MHM(Y,A)^{\rm nc}$ by external product which follows from [Ka2, Th.~4.4.1].
This finishes the proof of Theorem~2.
\ms\nin
{\bf 4.10.~Proof of Corollary~1.} This follows from the construction of the standard functors in [Sa3]. (For instance, we use Beilinson-type resolutions [Be1] for the definition of direct images.) The details are left to the reader.
\ms\nin
{\bf 4.11.~Remarks.} (i) The definition (or characterization) of mixed Hodge modules in this paper has been used in some papers (see e.g. Inv.\ Math.\ 144 (2001), p.~548). In fact, a similar argument is needed in order to show the stability by external products, and it seems simpler to take this as the definition of mixed Hodge modules from the beginning.
\ms
(ii) The theory of mixed Hodge modules in the analytic case is quite different from the algebraic case. For instance, a mixed Hodge module whose underlying perverse sheaf is a locally constant sheaf on a complex manifold is just a graded-polarizable variation of mixed Hodge structure (unless some extendability condition is imposed).
We have only cohomological functors $\Hc^jf_*$, etc.\ in general, and it is quite unclear whether the derived categories of mixed Hodge modules are really good categories.
\sk
In the analytic case it seems possible to give a definition similar to the one in this paper by assuming that variations of mixed structure are admissible with respect to the partial compactification given by the support of weakly mixed Hodge modules.
\ms
(iii) In [Sa3], the category of mixed Hodge modules was defined by restricting the category infinitely many times by using the stability by the nearby and vanishing cycle functors and the open direct images for the complements of locally principal divisors, and moreover by the pull-backs by smooth projections. However, the last one was unnecessary, since we have already restricted the category more than sufficiently and it was actually imposed only for a trivial proof of the stability by this functor. Moreover it is effectively enough to take the restrictions only $n$ times where $n$ is the dimension of the support.
Indeed, these imply first that a module is generically an admissible variation of mixed Hodge structure (see Remark~(iv) below).
To see that the whole module is a mixed Hodge module, it is then enough to verify the well-definedness along some function locally in order to proceed by induction on the dimension of the support (using the new definition in this paper).
\ms
(iv) The proof of the above assertion is not completely trivial. For instance, if the variety is non-normal, then we have to show the following: Let $\M$ be a weakly mixed Hodge module on a curve $X$, which is a variation of mixed Hodge structure on the complement $U$ of a singular point of $X$. Let $g$ be a function on $X$ vanishing only at this point. Assume the nearby and vanishing cycle functors along $g=0$ are well-defined for $\M$. Then $\M|_U$ is admissible at this point. (This can be shown by using the direct image by the function.)

\end{document}